\documentclass[preprint,11pt,numbers,nopreprintline]{elsarticle}
\usepackage[a4paper,margin=25mm]{geometry}
\usepackage{lmodern}
\usepackage{amsmath,amssymb,amsthm,mathtools,microtype}
\usepackage[hidelinks]{hyperref}
\hypersetup{pdftitle={When is the deepest drawdown also the longest?}, pdfauthor={Yinhao Wu}}
\journal{}

\newtheorem{theorem}{Theorem}
\newtheorem{lemma}{Lemma}
\newtheorem{proposition}{Proposition}
\newtheorem{corollary}{Corollary}
\theoremstyle{remark}

\newcommand{\Pp}{\mathbb P}
\newcommand{\Ee}{\mathbb E}
\newcommand{\Nn}{\mathbb N}
\newcommand{\ind}{\mathbf 1}
\newcommand{\asto}{\xrightarrow{\mathrm{a.s.}}}
\newcommand{\pto}{\xrightarrow{\Pp}}
\newcommand{\dto}{\xrightarrow{d}}
\newcommand{\eps}{\varepsilon}
\newcommand{\ccyc}{c_{\mathrm{cyc}}}

\begin{document}
\begin{frontmatter}
\title{When is the deepest drawdown also the longest?}
\author{Yinhao Wu}
\ead{yinhaowu@stu.sufe.edu.cn}
\affiliation{organization={School of Mathematics, Shanghai University of Finance and Economics}, city={Shanghai}, postcode={200433}, country={China}}

\begin{abstract}
Does the deepest drawdown take the longest to recover? For random walks with iid increments of positive finite mean, we find that the deepest and longest episodes are eventually distinct almost surely under light-tailed losses, but coincide with probability tending to one under heavy-tailed losses. More precisely, when losses have a finite exponential moment, the deepest episode's duration, relative to the longest among $n$ episodes, converges in probability to an explicit $\rho\in(0,1)$, while its descending duration rank is $n^{1-\rho+o_{\Pp}(1)}$. These limits hold for every fixed number of leading depths; in the Gaussian model, $\rho=1/2$, independently of drift and variance. The separation reflects the different exponential costs of reaching a large depth and sustaining a long episode. For regularly varying losses, a common large loss governs both depth and duration: the duration ratio tends to one almost surely, and for every fixed $k$, the $k$ deepest and $k$ longest episodes have the same indices in the same order with probability tending to one. For any $k_n\to\infty$ with $k_n=o(n)$, the proportion of common episodes among the $k_n$ deepest and $k_n$ longest converges in probability to zero in the light-tailed case and to one in the heavy-tailed case. These ranking laws also hold for completed drawdowns within a deterministic observation horizon.
\end{abstract}

\end{frontmatter}

\section{Introduction}

An investment can grow over the long run while spending extended periods below a previous peak. Both the size of the loss and the wait for recovery matter to investors: drawdown depth records the largest decline from the peak during an episode, while duration measures the time from the peak to its first recovery. A deeper loss leaves a larger deficit to recover, suggesting that the deepest drawdowns should also be among the longest. Yet a succession of modest setbacks can prolong an episode without producing an exceptionally deep loss. Mahmoud~\cite{Mahmoud2017} distinguishes the longest drawdown duration from the duration of the deepest drawdown in the study of temporal risk.

We study when the deepest drawdown is also the longest for a random walk with iid increments of positive finite mean, which can represent cumulative log returns. As the number of completed episodes grows, the answer depends on the loss tail: when losses have a finite exponential moment, the deepest and longest episodes are eventually distinct almost surely; for regularly varying losses with index greater than one, they coincide with probability tending to one. We quantify this contrast through the deepest episode's duration rank and the overlap between lists of leading episodes.

The magnitude of drawdowns has been studied through both fixed-horizon maxima and passage times at prescribed levels. For Brownian motion with drift, Magdon-Ismail et al.~\cite{MagdonIsmailEtAl2004} derive the distribution and expected value of the maximum drawdown. In positive-drift random walks, Meilijson~\cite{Meilijson2009} bounds the expected running maximum before a prescribed drawdown in terms of the adjustment coefficient. Pospisil, Vecer and Hadjiliadis~\cite{PospisilVecerHadjiliadis2009} obtain joint laws of the running extrema at a drawdown passage time for general diffusions, and compute the probability that a prescribed drawdown precedes a prescribed drawup. Joint passage laws also describe drawdowns of spectrally negative L\'evy processes~\cite{MijatovicPistorius2012}, with an integral-equation approach available for time-homogeneous Markov processes~\cite{LandriaultLiZhang2017Markov}. For future drawdown extremes of L\'evy processes, Baurdoux, Palmowski and Pistorius~\cite{BaurdouxPalmowskiPistorius2017} derive tail asymptotics under Cram\'er and convolution-equivalent assumptions.

The time spent below a peak is another measure of drawdown risk. The time from the last peak to a drawdown of a given size has been studied under diffusion dynamics~\cite{ZhangHadjiliadis2012} and through path decompositions of spectrally negative L\'evy processes~\cite{VardarAcarCaglarAvram2021}. Landriault, Li and Zhang~\cite{LandriaultLiZhang2017Duration} derive the Laplace transform of the first time a drawdown has lasted a prescribed period for a class of L\'evy models. Their study of repeated Brownian drawdowns~\cite{LandriaultLiZhang2015} gives frequency laws with and without recovery between successive drawdown times.

Regeneration at recovery makes successive completed drawdowns independent and identically distributed. Depth and duration remain dependent within each episode, and their joint extreme behavior determines how the two rankings are related. To quantify this relation, let $L_i$ be the duration of episode $i$ and $I_n$ the index of the deepest among the first $n$ completed nonzero episodes. Write $R_n$ for the descending duration rank of that episode. Thus $R_n=1$ means that the deepest episode is also the longest, while a growing rank records its separation from the longest episodes.

For light-tailed losses, both the largest depth and the longest duration grow logarithmically with the number of episodes. Reaching a large depth and sustaining a long episode have different exponential costs, and this difference separates their rankings. Assuming a finite exponential moment, we show that an explicit constant $\rho\in(0,1)$ satisfies
\begin{equation*}
 \frac{L_{I_n}}{\max_{i\le n}L_i}\pto\rho,
 \qquad R_n=n^{1-\rho+o_{\Pp}(1)}.
\end{equation*}
The duration ratio quantifies the difference between the two extremes, while the rank exponent describes the growth rate of the number of episodes ranked ahead of the deepest. For every fixed $k$, the same limits hold jointly for the $k$ deepest episodes, and these episodes eventually lie outside any fixed list of the longest, almost surely. In the Gaussian model, $\rho=1/2$ for every positive mean and variance: the deepest drawdown lasts asymptotically half as long as the longest, and its duration rank has logarithmic growth exponent $1/2$.

Under regular variation, a single large negative increment supplies both the depth of an extreme episode and the deficit that the subsequent positive drift must recover. The deepest episode lasts asymptotically as long as the longest, and the two episodes coincide with probability tending to one:
\begin{equation*}
 \frac{L_{I_n}}{\max_{i\le n}L_i}\asto1,
 \qquad
 \Pp(R_n=1)\longrightarrow1.
\end{equation*}
More generally, for every fixed $k$, the $k$ deepest and $k$ longest episodes have the same indices in the same order with probability tending to one. On the scale of the largest depth, the difference between depth and drift-adjusted duration vanishes uniformly across the sample, almost surely.

When the number of selected episodes grows, the same contrast appears in the overlap between the two leading lists. For every $k_n\to\infty$ with $k_n=o(n)$, the proportion of common episodes among the $k_n$ deepest and $k_n$ longest converges in probability to zero for light-tailed losses and to one under regular variation.

The excursion representation links drawdown extremes to Iglehart's queueing analysis~\cite{Iglehart1972}. For the simple symmetric random walk, exact and limiting laws describe ranked excursion heights and lengths~\cite{CsakiHu2003}. Perfilev and Wachtel~\cite{PerfilevWachtel2022} derive joint local asymptotics for an excursion's maximum, its time of attainment, and its length, including the length conditional on a large maximum.

Under light tails, we first determine the duration scale of episodes selected by depth, using a decomposition at the deepest point to control the subsequent recovery. Counting episodes with longer durations then gives the rank exponent, and a joint exponential bound yields eventual separation almost surely. Under regular variation, a weighted bound on depth--duration deviations yields an almost-sure uniform approximation across the sample. Comparing this error with the gaps between the largest depths identifies the same episode indices in the same order, with probability tending to one. For growing lists, we use the same estimates to compare joint and marginal exceedance counts and determine the overlap proportion.

Section~\ref{sec:main} defines the drawdowns and their ranks and gives the common excursion representation. Section~\ref{sec:results} presents the two ranking laws and the overlap limit for growing leading lists, followed by examples. Sections~\ref{sec:light} and~\ref{sec:heavy} give the proofs for light- and heavy-tailed losses, respectively. Section~\ref{sec:calendar} transfers both ranking laws to completed drawdowns within a deterministic observation horizon.

\section{Model and drawdown ranks}\label{sec:main}

\subsection{Completed episodes and their ranks}
Let $X_1,X_2,\ldots$ be iid real random variables such that
\begin{equation}\label{eq:common}
 \Ee|X_1|<\infty,\qquad \mu:=\Ee X_1>0,\qquad
 q:=\Pp(X_1<0)>0.
\end{equation}
Set $S_0=0$, $S_n=\sum_{j=1}^nX_j$, and $M_n=\max_{0\le j\le n}S_j$. The drawdown is
\begin{equation*}
 D_n=M_n-S_n=(D_{n-1}-X_n)^+.
\end{equation*}
When $S_n$ represents cumulative log returns, percentage drawdown is $1-e^{-D_n}$; this increasing transformation preserves the depth rankings and all recovery times.

We count nonzero completed episodes. Starting from $r_0=0$, define
\begin{align*}
 g_i&=\inf\{n>r_{i-1}:X_n<0\}-1,&
 r_i&=\inf\{n>g_i:S_n\ge S_{g_i}\},\\
 L_i&=r_i-g_i,&
 H_i&=S_{g_i}-\min_{g_i\le j\le r_i}S_j.
\end{align*}
Between $r_{i-1}$ and $g_i$ the path stays at its current running maximum; $g_i$ is the last peak immediately preceding a loss. The endpoint $r_i$ is the first recovery of that peak, with equality counted as recovery. Thus $H_i>0$ and $L_i\ge2$. The strong law ensures that each recovery is finite. The duration includes the descent to the trough as well as the subsequent recovery.

For $1\le j\le n$, let $I_{n,j}$ be the index of the $j$th largest $H_i$ among $i\le n$, and let $J_{n,j}$ be the index of the $j$th largest $L_i$. In each ordering, ties are broken by the smaller index. Define
\begin{equation*}
 I_n=I_{n,1},\qquad H_n^*=\max_{i\le n}H_i,\qquad
 L_n^*=\max_{i\le n}L_i.
\end{equation*}
The duration rank of the $j$th deepest episode is
\begin{equation*}
 R_{n,j}=1+\sum_{i=1}^n\ind_{\{L_i>L_{I_{n,j}}
       \ \text{or}\ (L_i=L_{I_{n,j}},\ i<I_{n,j})\}},
 \qquad R_n=R_{n,1}.
\end{equation*}
Thus $R_n=1$ exactly when the deepest episode is ranked first by duration, that is, $I_n=J_{n,1}$. For $1\le k\le n$, the fraction of the $k$ deepest episodes that also belong to the $k$ longest is
\[
 O_{n,k}=\frac1k\sum_{j=1}^k\ind_{\{R_{n,j}\le k\}}.
\]

\subsection{Excursion representation and regeneration}\label{sec:regeneration}

Viewed in the direction of losses, a drawdown is a positive excursion of a walk with drift $-\mu$, ending when the loss walk first returns to or falls below zero. Write $Y_j=-X_j$, $Z_n=\sum_{j=1}^nY_j$, and
\begin{equation}\label{eq:aux}
 \tau=\inf\{n\ge1:Z_n\le0\},\qquad
 H^0=\max_{0\le j<\tau}Z_j.
\end{equation}
Here $H^0$ and $\tau$ are the height and length of the unconditioned excursion. A nonzero episode starts with an increment having the law of $X_1$ given $X_1<0$, followed by increments with the original iid law. Thus $(H_i,L_i)$ has the law of $(H^0,\tau)$ conditional on $Y_1>0$. Restarting at the recovery stopping times $r_i$ makes the successive pairs independent and identically distributed.

To relate episode counts to observation time, write $W_i=g_i-r_{i-1}$ for the wait at the running maximum before episode $i$. The triples $(W_i,H_i,L_i)$ are iid, and $W_i$ is independent of $(H_i,L_i)$, with $\Pp(W_i=j)=(1-q)^jq$ for $j\ge0$. The mean inter-completion time is
\begin{equation}\label{eq:mean-cycle}
 \ccyc:=\Ee(W_1+L_1)=\frac{\Ee\tau}{q}<\infty.
\end{equation}
To verify finiteness under \eqref{eq:common}, choose a truncation level $B$ large enough that $\Ee\max(Y_1,-B)=-\mu_B<0$. Let $\tau_B$ be the first weak exit from $(0,\infty)$ of the walk with these truncated increments. Pathwise $\tau\le\tau_B$. The truncated walk is at least $-B$ at $\tau_B\wedge n$, since a downward crossing overshoots by at most $B$. Wald's identity at this bounded stopping time yields $-\mu_B\Ee(\tau_B\wedge n)\ge-B$. Monotone convergence gives $\Ee\tau\le B/\mu_B$. Finally, $\Ee\tau=(1-q)+q\Ee L_1$ and $\Ee W_1=(1-q)/q$, proving \eqref{eq:mean-cycle}.

\section{Two ranking laws}\label{sec:results}

Under the common assumptions \eqref{eq:common}, we study two cases for the loss tail.

\subsection{Light-tailed losses}\label{sec:light-results}
Assume that the loss distribution has a finite exponential moment:
\begin{equation*}\tag{L}\label{ass:light}
 \Ee e^{-\eta X_1}<\infty\qquad\text{for some }\eta>0.
\end{equation*}
Define
\begin{equation*}
 \begin{gathered}
 \psi(\theta)=\log\Ee e^{-\theta X_1},\qquad \theta\ge0,\\[3pt]
 \alpha=\sup\{\theta\ge0:\psi(\theta)\le0\},\qquad
 \beta=-\min_{0\le\theta\le\alpha}\psi(\theta).
 \end{gathered}
\end{equation*}
The parameters $\alpha$ and $\beta$ are the exponential decay rates of the depth and duration tails, respectively. Both are positive and finite, and $\psi(\alpha)\le0$. Set
\begin{equation*}
 v=\begin{cases}
     \Ee[-X_1e^{-\alpha X_1}],&\psi(\alpha)=0,\\
     \infty,&\psi(\alpha)<0,
   \end{cases}
 \qquad A=\frac1\mu+\frac1v,
\end{equation*}
where the expectation may be infinite and $1/\infty=0$. When $v<\infty$, $\alpha$ is the Cram\'er root and $v$ is the mean loss under the exponential tilt with density $e^{-\alpha X_1}$. A deep episode typically descends at speed $v$ and recovers at drift $\mu$, so an episode of depth $h$ lasts approximately $Ah$ periods. When $v=\infty$, the descent takes $o(h)$ periods and recovery determines the first-order duration $h/\mu$.

\begin{theorem}[Rank separation]\label{thm:light}
Under \eqref{eq:common} and \eqref{ass:light}, define
\begin{equation*}
 \rho=\frac{\beta A}{\alpha}
     =\frac{\beta}{\alpha}\left(\frac1\mu+\frac1v\right).
\end{equation*}
Then $0<\rho<1$, and
\begin{equation}\label{eq:rank-limit}
 \frac{L_{I_n}}{L_n^*}\pto\rho,
 \qquad
 \frac{\log R_n}{\log n}\pto1-\rho.
\end{equation}
Almost surely, $\limsup_{n\to\infty}L_{I_n}/L_n^*=\rho$, and the deepest and longest episodes are different for all sufficiently large $n$. For every fixed positive integer $k$, the limits in \eqref{eq:rank-limit} hold jointly with $(I_n,R_n)$ replaced by $(I_{n,j},R_{n,j})$, $1\le j\le k$. The index sets of the $k$ deepest and $k$ longest episodes are eventually disjoint almost surely.

The scales underlying these limits are
\begin{equation}\label{eq:light-max}
 \frac{H_n^*}{\log n}\asto\frac1\alpha,
 \qquad
 \frac{L_n^*}{\log n}\asto\frac1\beta,
 \qquad
 \frac{L_{I_n}}{\log n}\pto\frac A\alpha.
\end{equation}
\end{theorem}

A deep episode has exponential cost $\alpha/A$ per unit duration, whereas a long episode has cost $\beta$. Strict convexity gives $\beta<\alpha/A$, and hence $\rho<1$. The duration-tail rate in Lemma~\ref{lem:light-tails} also gives
\begin{equation*}
 n\,\Pp\!\left(L_1>\frac A\alpha\log n\right)
   =n^{1-\rho+o(1)}.
\end{equation*}
This expected count at the selected duration scale explains the rank exponent.

\subsection{Heavy-tailed losses}\label{sec:heavy-results}
Assume that the loss tail is regularly varying:
\begin{equation*}\tag{H}\label{ass:heavy}
 \overline F(x):=\Pp(X_1<-x)=x^{-\gamma}\ell(x),
 \qquad \gamma>1,
\end{equation*}
where $\ell$ is slowly varying, meaning that $\ell(tx)/\ell(x)\to1$ for each fixed $t>0$.

Here a single large loss creates both the depth of an extreme episode and the deficit that the positive drift must recover. The resulting relation between depth and duration also determines which episodes occupy the leading ranks.

\begin{theorem}[Rank agreement]\label{thm:heavy}
Under \eqref{eq:common} and \eqref{ass:heavy}, for every fixed positive integer $k$,
\begin{equation}\label{eq:topk}
 \Pp(I_{n,j}=J_{n,j}\ \text{for all }1\le j\le k)
 \longrightarrow1.
\end{equation}
In particular,
\begin{equation}\label{eq:heavy-conseq}
 \Pp(R_n=1)\longrightarrow1,\qquad
 \frac{L_{I_n}}{L_n^*}\asto1.
\end{equation}

Across the sample, depth and drift-adjusted duration satisfy
\begin{equation}\label{eq:uniform-alignment}
 \frac1{H_n^*}\max_{i\le n}|H_i-\mu L_i|\asto0.
\end{equation}

Moreover, choose $b_n\to\infty$ such that
\begin{equation*}
 n\ccyc\overline F(b_n)\longrightarrow1.
\end{equation*}
Let $\Gamma_1<\Gamma_2<\cdots$ be the successive arrival times of a Poisson process with rate one. The leading depths and drift-adjusted durations have the joint limit
\begin{equation}\label{eq:order-limit}
 \left(\frac{H_{I_{n,j}}}{b_n},
       \frac{\mu L_{I_{n,j}}}{b_n}\right)_{j=1}^{k}
 \dto\left(\Gamma_j^{-1/\gamma},\Gamma_j^{-1/\gamma}\right)_{j=1}^{k}.
\end{equation}
\end{theorem}

For each fixed $k$, the gaps between the largest depths, including the one between the $k$th and $(k+1)$st, dominate the uniform error in \eqref{eq:uniform-alignment} with probability tending to one. Replacing depth by $\mu$ times duration therefore preserves the leading episode indices and their order.

\subsection{Growing leading lists}\label{sec:overlap}

We turn from a fixed number of leading episodes to lists whose size $k_n$ grows with the sample. The overlap $O_{n,k_n}$ measures the fraction of the $k_n$ deepest episodes that also belong to the $k_n$ longest.

\begin{corollary}[Overlap of growing lists]\label{cor:overlap}
Under \eqref{eq:common}, let $k_n$ be deterministic positive integers such that $k_n\to\infty$ and $k_n/n\to0$. Then
\begin{equation}\label{eq:overlap}
 O_{n,k_n}\pto
 \begin{cases}
  0,&\text{under \eqref{ass:light}},\\
  1,&\text{under \eqref{ass:heavy}}.
 \end{cases}
\end{equation}
If, in addition, $k_n/\log n\to\infty$, both limits hold almost surely.
\end{corollary}

\subsection{Examples}\label{sec:examples}

We illustrate the two ranking laws with three examples. For the Gaussian and two-point models, the limiting duration ratio has an explicit form; the shifted Pareto model falls under the heavy-tailed case.

\paragraph{Example 1 (Gaussian model).}
Let $X_1\sim N(\mu,\sigma^2)$ with $\mu,\sigma>0$. Then
\[
 \begin{gathered}
 \psi(\theta)=-\mu\theta+\tfrac12\sigma^2\theta^2,\\[3pt]
 \alpha=\frac{2\mu}{\sigma^2},\quad
 v=\mu,\quad
 \beta=\frac{\mu^2}{2\sigma^2}.
 \end{gathered}
\]
Therefore $\rho=1/2$, and Theorem~\ref{thm:light} gives
\begin{equation*}
 \frac{L_{I_n}}{L_n^*}\pto\frac12,
 \qquad R_n=n^{1/2+o_{\Pp}(1)}.
\end{equation*}
The absolute scales depend on drift and volatility, while the duration ratio and rank exponent do not.

The factor $1/2$ also has a direct cost interpretation. Since $v=\mu$, a deep descent and its recovery take about $2h/\mu$ periods in total, with exponential cost $\alpha h=2\mu h/\sigma^2$. A long excursion sustained by near-zero average motion incurs cost $\beta=\mu^2/(2\sigma^2)$ per period. At the same exponential cost, its duration is about $4h/\mu$, twice that of the deep episode.

\paragraph{Example 2 (Two-point model).}
For comparison, let $\Pp(X_1=1)=p>1/2$ and $\Pp(X_1=-1)=1-p$. Then
\[
 \alpha=\log\frac p{1-p},\qquad
 \mu=v=2p-1,\qquad
 \beta=-\log\bigl(2\sqrt{p(1-p)}\bigr).
\]
It follows that
\begin{equation*}
 \rho=
 \frac{-2\log\bigl(2\sqrt{p(1-p)}\bigr)}
 {(2p-1)\log[p/(1-p)]}.
\end{equation*}
For $p=0.7$, this is approximately $0.5144$, so the limiting duration ratio varies with the increment law.

\paragraph{Example 3 (Shifted Pareto model).}
For the heavy-tailed case, let
\[
 X_1=c-Z,\qquad \Pp(Z>z)=z^{-\gamma}\quad(z\ge1),
 \qquad \gamma>1,\quad c>\frac{\gamma}{\gamma-1}.
\]
Then
\[
 \mu=c-\frac{\gamma}{\gamma-1}>0,\qquad
 \overline F(x)=(x+c)^{-\gamma}\quad(x\ge0).
\]
The increments are integrable and bounded above, so Theorem~\ref{thm:heavy} applies. For $\gamma=3/2$ and $c=4$, the mean is $\mu=1$ and the loss has infinite second moment.

\section{Proofs for light-tailed losses}\label{sec:light}

Throughout this section, assume \eqref{eq:common} and \eqref{ass:light}. We first obtain the two marginal tail rates, then decompose the loss excursion at its maximum to study recovery after selection by depth. A joint exponential bound gives almost-sure separation and the overlap limit for growing lists.

\subsection{Depth and duration estimates}

The definition of $\alpha$ gives $0<\alpha<\infty$ and $\Ee e^{-\alpha X_1}\le1$: positivity follows from $\psi'(0+)=-\mu<0$, finiteness from $\Pp(X_1<0)>0$, and the endpoint inequality from Fatou's lemma. Dominated convergence makes $\psi$ continuous on $[0,\alpha]$. It is strictly negative on $(0,\alpha)$, so $\beta>0$.

For the loss walk $Z$ in \eqref{eq:aux}, write $T_h=\inf\{n\ge1:Z_n\ge h\}$. For $0\le\theta\le\alpha$, let $\Pp_\theta$ give the loss increments the iid tilted law
\[
 \Pp_\theta(Y_1\in dy)
   =e^{\theta y-\psi(\theta)}\Pp(Y_1\in dy).
\]
When $v<\infty$, $\psi(\alpha)=0$ and this law has mean $v>0$ at $\theta=\alpha$.

\begin{lemma}\label{lem:light-tails}
Under \eqref{eq:common} and \eqref{ass:light},
\begin{equation}\label{eq:light-tails}
 \lim_{h\to\infty}\frac1h\log\Pp(H_1\ge h)=-\alpha,
 \qquad
 \lim_{t\to\infty}\frac1t\log\Pp(L_1>t)=-\beta.
\end{equation}
For all $h>0$, $\Pp(H_1\ge h)\le q^{-1}e^{-\alpha h}$. If $v<\infty$, there are constants $c,C>0$ such that, for all sufficiently large $h$,
\begin{equation}\label{eq:depth-bounds}
 ce^{-\alpha h}\le\Pp(H_1\ge h)\le Ce^{-\alpha h}.
\end{equation}
\end{lemma}
\begin{proof}
We use $(H^0,\tau)$ from \eqref{eq:aux}; conditioning on $Y_1>0$ divides positive-height and long-duration probabilities by $q$. Let $\theta_*$ minimize $\psi$ on $[0,\alpha]$. The nonnegative supermartingale $e^{\alpha Z_n}$ and Chernoff's inequality give
\[
 \Pp(H^0\ge h)\le e^{-\alpha h},\qquad
 \Pp(\tau>n)\le\Pp(Z_n>0)\le e^{-\beta n}.
\]

Suppose first that $v<\infty$. Summing the change of measure over the possible values of $T_h$ yields
\begin{equation}\label{eq:tilt-stopped}
 \Pp(H^0\ge h)=\Pp(T_h<\tau)
 =\Ee_\alpha[e^{-\alpha Z_{T_h}};\ T_h<\tau].
\end{equation}
For $p(x)=\Pp(\sup_{n\ge0}Z_n>x)$, the classical Cram\'er estimate gives $e^{\alpha x}p(x)\to c_0\in(0,\infty)$, along the support lattice in the arithmetic case. Decomposition at the first weak exit, as in Iglehart~\cite{Iglehart1972}, gives
\[
 p(x)=\Pp(H^0>x)+\Ee[\ind_{\{H^0\le x\}}p(x-Z_\tau)].
\]
Since $Z_\tau\le0$ and $p(u)\le e^{-\alpha u}$, dominated convergence yields
\[
 e^{\alpha x}\Pp(H^0>x)
 \longrightarrow c_0\bigl(1-\Ee e^{\alpha Z_\tau}\bigr)>0.
\]
The strict inequality follows from $\Pp(Y_1<0)>0$. Rounding to adjacent lattice points when needed proves \eqref{eq:depth-bounds}.

For the duration lower bound, use $\theta_*\in(0,\alpha)$: the tilted increments under $\Pp_{\theta_*}$ have mean zero and finite positive variance. Choose $0<a_0<a_1$. A vector with positive total sum has a cyclic shift with strictly positive partial sums: start after the last minimum of its partial sums. Cyclic invariance and the central limit theorem therefore give
\[
 \begin{aligned}
 &\Pp_{\theta_*}(\tau>n,\ a_0\sqrt n<Z_n<a_1\sqrt n)\\
 &\qquad\ge\frac1n\Pp_{\theta_*}(a_0\sqrt n<Z_n<a_1\sqrt n)
 \ge\frac{c}{n}.
 \end{aligned}
\]
Changing measure at time $n$ yields
\[
 \Pp(\tau>n)\ge\frac cn e^{-\beta n-\theta_*a_1\sqrt n},
\]
which proves the duration rate for $v<\infty$.

For the general case, truncate the loss increments to $\min(Y_j,b)$, $b>0$. The truncated walk has negative mean, a finite Cram\'er root $\alpha_b$, and a finite tilted mean there. Its excursion height and duration are bounded by $H^0$ and $\tau$, respectively. Write
\[
 \phi_b(\theta)=\Ee e^{\theta\min(Y_1,b)},\qquad
 e^{-\beta_b}=\min_{\theta\ge0}\phi_b(\theta).
\]
The roots $\alpha_b$ decrease to $\alpha$: if their limit were larger, any $\theta$ strictly between the two would satisfy $\phi_b(\theta)\le1$ for every $b$, hence $\Ee e^{\theta Y_1}\le1$ by monotone convergence.

To identify the limiting minima, choose $d>0$ with $\Pp(Y_1\ge d)>0$. For $b\ge d$, the bound $\phi_b(\theta)\ge\Pp(Y_1\ge d)e^{\theta d}$ confines all minimizers to a fixed compact interval. If a subsequence of minimizers $\theta_b$ converges to $\theta_0$, Fatou's lemma and $\phi_b\le e^\psi$ give
\[
 e^{-\beta}
 \le\Ee e^{\theta_0Y_1}
 \le\liminf_b\phi_b(\theta_b)
 \le\limsup_b\phi_b(\theta_b)
 \le e^{-\beta}.
\]
Thus $\beta_b\downarrow\beta$. The finite-$v$ result applied to the truncated excursions now gives
\[
 \liminf_{h\to\infty}\frac1h\log\Pp(H^0\ge h)\ge-\alpha_b,
 \qquad
 \liminf_{n\to\infty}\frac1n\log\Pp(\tau>n)\ge-\beta_b.
\]
Letting $b\to\infty$ and combining with the upper bounds proves \eqref{eq:light-tails}.
\end{proof}

We combine these marginal rates with joint estimates to determine the duration of the deepest episodes.

\begin{proposition}\label{prop:joint-light}
Under \eqref{eq:common} and \eqref{ass:light}, for each $\eps>0$ there are constants $C_\eps,c_\eps>0$ such that
\begin{equation}\label{eq:joint-long}
 \Pp(H_1\ge h,\ L_1>(A+\eps)h)
 \le C_\eps e^{-(\alpha+c_\eps)h},\qquad h\ge1.
\end{equation}
There is also a bounded function $r_\eps(x)\to0$ as $x\to\infty$ such that
\begin{equation}\label{eq:recovery-lower}
 \Pp\bigl(L_1<(\mu^{-1}-\eps)H_1\mid H_1\bigr)
 \le r_\eps(H_1)\quad\text{almost surely}.
\end{equation}
If $v<\infty$, then
\begin{equation}\label{eq:joint-light}
 \Pp(H_1\ge h,\ |L_1-Ah|>\eps h)=o(e^{-\alpha h}),
 \qquad h\to\infty.
\end{equation}
\end{proposition}
\begin{proof}
Work with $(H^0,\tau)$ and let $G$ be the first time the loss excursion attains its maximum. Positive-height probabilities are divided by $q$ when passing to a nonzero episode. For the original walk $S$, put $\sigma_x=\inf\{n\ge1:S_n\ge x\}$ and
\[
 Q_x(B)=\Pp(\sigma_x\in B,\ S_j\ge0\text{ for }j<\sigma_x),
 \qquad x>0,\quad B\subseteq\Nn.
\]
For each fixed $k\ge1$, define
\[
 \nu_k(dx)=\Pp(\tau>k,\ \max_{0\le j<k}Z_j<Z_k,\ Z_k\in dx).
\]
Independence after time $k$ gives
\[
 \Pp(H^0\in dx,G=k,\tau-G\in B)=\nu_k(dx)Q_x(B).
\]
The strict inequality before $k$ and the weak inequality in $Q_x$ encode first attainment of the maximum. Summing over $k$ identifies the conditional law of $\tau-G$ given $H^0=x$ as $Q_x/Q_x(\Nn)$.

Let $p_*=\Pp(S_j\ge0\text{ for all }j\ge1)$. Choose $a>0$ with $\Pp(X_1\ge a)>0$. Since $\inf_{j\ge0}S_j> -\infty$ almost surely, independence of the increments gives, for sufficiently large $m$,
\[
 p_*\ge\Pp(X_1\ge a)^m\Pp(\inf_{j\ge0}S_j\ge-ma)>0.
\]
Hence $Q_x(\Nn)\ge p_*$ and, for almost every $x>0$ under the height law,
\[
 \Pp(\tau-G\in B\mid H^0=x)
 =\frac{Q_x(B)}{Q_x(\Nn)}
 \le p_*^{-1}\Pp(\sigma_x\in B).
\]
This kernel is unchanged on conditioning the episode to be nonzero. The strong law gives $\sigma_x/x\to1/\mu$ almost surely, so \eqref{eq:recovery-lower} holds with $r_\eps(x)=p_*^{-1}\Pp(\sigma_x<(\mu^{-1}-\eps)x)$. For $t>1/\mu$, Chernoff's bound gives, for sufficiently small $s>0$,
\[
 \Pp(\sigma_x>tx)
 \le e^{sx+\lfloor tx\rfloor\psi(s)}\le Ce^{-cx},
\]
since $\psi(s)=-\mu s+o(s)$. In particular, for each $\delta>0$,
\[
 \begin{aligned}
 &\Pp(H^0\ge h,\ \tau-G>(\mu^{-1}+\delta)H^0)\\
 &\quad\le p_*^{-1}\int_{[h,\infty)}
   \Pp(\sigma_x>(\mu^{-1}+\delta)x)\,\Pp(H^0\in dx)\\
 &\quad\le Ce^{-ch}\Pp(H^0\ge h)
 \le Ce^{-(\alpha+c)h}.
 \end{aligned}
\]

To bound the time $G$ spent reaching the maximum, use the identity
\[
 \lim_{\theta\uparrow\alpha}
 \frac{\alpha-\theta}{-\psi(\theta)}=\frac1v.
\]
Indeed, if $\psi(\alpha)=0$, use $\psi'(\alpha-)=v$, allowing an infinite left derivative; if $\psi(\alpha)<0$, the numerator tends to zero and the denominator to a positive value. For $t>1/v$, choose $\theta\in(0,\alpha)$ with $\theta-t\psi(\theta)>\alpha$. Then
\[
 \Pp(H^0\ge h,G>th)
 \le\sum_{j>th}\Pp(Z_j\ge h)
 \le C_\theta e^{-[\theta-t\psi(\theta)]h}.
\]
Choose $\delta>0$ so that $(v^{-1}+\delta)+(\mu^{-1}+\delta)(1+\delta)<A+\eps$. On $h\le H^0\le(1+\delta)h$, the bounds $G\le(v^{-1}+\delta)h$ and $\tau-G\le(\mu^{-1}+\delta)H^0$ imply $\tau\le(A+\eps)h$. The preceding estimates therefore yield
\[
 \Pp(H^0\ge h,\tau>(A+\eps)h)
 \le e^{-\alpha(1+\delta)h}
      +C_1e^{-(\alpha+c_1)h}+C_2e^{-(\alpha+c_2)h},
\]
proving \eqref{eq:joint-long}.

For \eqref{eq:joint-light}, suppose $v<\infty$ and choose $\delta>0$ with $\delta/\mu<\eps/3$. Under $\Pp_\alpha$ the loss walk has finite positive mean $v$, so $T_h/h\to1/v$ almost surely. Equation~\eqref{eq:tilt-stopped} gives
\[
 \Pp(T_h<\tau,\ |T_h-h/v|>\eps h/3)
 \le e^{-\alpha h}\Pp_\alpha(|T_h-h/v|>\eps h/3)
 =o(e^{-\alpha h}),
\]
and
\[
 \Pp(T_h<\tau,\ Z_{T_h}>(1+\delta)h)
 \le e^{-\alpha(1+\delta)h}.
\]
The increments after $T_h$ are independent of the past. Conditional on $z=Z_{T_h}$, the remaining duration has law $\sigma_z$, and its strong law yields
\[
 \sup_{h\le z\le(1+\delta)h}
 \Pp(|\sigma_z-z/\mu|>\eps h/3)\longrightarrow0.
\]
Multiplying by $\Pp(T_h<\tau)\le e^{-\alpha h}$ bounds the corresponding joint event by $o(e^{-\alpha h})$. Outside these three events,
\[
 |\tau-Ah|
 \le |T_h-h/v|+|\tau-T_h-z/\mu|+(z-h)/\mu\le\eps h.
\]
This proves \eqref{eq:joint-light}.
\end{proof}

\subsection{Ranks of the deepest episodes}\label{sec:light-ranks}

\begin{proof}[Proof of Theorem~\ref{thm:light}]
Let $\theta_*$ minimize $\psi$ on $[0,\alpha]$. If $v<\infty$, strict convexity gives
\[
 -\beta=\psi(\theta_*)>-\mu\theta_*,\qquad
 -\beta=\psi(\theta_*)>v(\theta_*-\alpha).
\]
Dividing by $\mu$ and $v$ and adding gives $\beta A<\alpha$. If $v=\infty$, strict Jensen gives $\beta<\mu\theta_*\le\mu\alpha$, again yielding $\beta A<\alpha$. Thus $0<\rho<1$ in both cases.

We use the marginal tail rates to control maxima and exceedance counts. For iid nonnegative $B_i$ with $\log\Pp(B_1>x)/x\to-t$, $t>0$, the maximum satisfies
\begin{equation}\label{eq:max-log}
 \frac{\max_{i\le n}B_i}{\log n}\asto\frac1t.
\end{equation}
To see this and obtain the counts needed below, fix $0<b<1/t$ and put $C_n(b)=\sum_{i\le n}\ind_{\{B_i>b\log n\}}$. Its mean is $n^{1-tb+o(1)}$, and Chernoff's inequality gives
\[
 \Pp\bigl(|C_n(b)-\Ee C_n(b)|>\tfrac12\Ee C_n(b)\bigr)
 \le2\exp[-\Ee C_n(b)/12].
\]
The right-hand side is summable over $n$, so
\begin{equation}\label{eq:counts}
 \frac{\log(1+C_n(b))}{\log n}\asto1-tb.
\end{equation}
This also holds for weak exceedances, since shifting the threshold by one preserves the exponential rate. In particular, the $k$th largest value exceeds $(1-\delta)\log n/t$ eventually, for every fixed $k$ and $\delta\in(0,1)$. For the upper bound,
\[
 \sum_{m\ge1}\Pp\!\left(
  \max_{i\le2^{m+1}}B_i>\frac{(1+\delta)m\log2}{t}\right)
 \le\sum_{m\ge1}2^{m+1}
       \Pp\!\left(B_1>\frac{(1+\delta)m\log2}{t}\right)<\infty.
\]
Borel--Cantelli on $2^m\le n\le2^{m+1}$, followed by $\delta\downarrow0$, proves \eqref{eq:max-log}, also for each fixed upper order statistic. Applying these facts with \eqref{eq:light-tails} gives
\[
 \frac{H_{I_{n,k}}}{\log n}\asto\frac1\alpha,
 \qquad \frac{L_{J_{n,k}}}{\log n}\asto\frac1\beta.
\]
The one-sided estimate \eqref{eq:joint-long} gives an almost-sure upper bound on the selected durations. Fix $\eps>0$ and choose $\delta\in(0,1)$ so that $(\alpha+c_\eps)(1-\delta)>\alpha$. Put $h_m=(1-\delta)m\log2/\alpha$. Then
\[
 \Pp\bigl(\exists i\le2^{m+1}:H_i\ge h_m,\ L_i>(A+\eps)h_m\bigr)
 \le C_\eps2^{m+1}e^{-(\alpha+c_\eps)h_m},
\]
which is summable in $m$. Among the first $2^m$ observations, the number of depths at least $h_m$ has mean $2^{\delta m+o(m)}$ and is at least $k$ eventually almost surely. Throughout $2^m\le n\le2^{m+1}$, all of the first $k$ selected depths therefore exceed $h_m$ and their durations are at most $(A+\eps)h_m$, eventually almost surely. Letting $\eps\downarrow0$ yields
\begin{equation}\label{eq:transfer-upper}
 \limsup_{n\to\infty}\max_{j\le k}\frac{L_{I_{n,j}}}{\log n}
 \le\frac A\alpha\quad\text{almost surely}.
\end{equation}

For convergence in probability when $v<\infty$, fix $\eps,M>0$ and put $h_n=(\log n)/\alpha-M$. Equation~\eqref{eq:joint-light} gives
\[
 \Pp\bigl(\exists i\le n:H_i\ge h_n,\ |L_i-Ah_n|>\eps h_n\bigr)
 \le n\,o(e^{-\alpha h_n})\longrightarrow0.
\]
The number of depths at least $h_n$ has mean at least $ce^{\alpha M}$ by \eqref{eq:depth-bounds}. For $M$ large enough that $ce^{\alpha M}\ge2k$, a binomial bound yields
\[
 \limsup_{n\to\infty}\Pp(H_{I_{n,k}}<h_n)
 \le\exp[-ce^{\alpha M}/8].
\]
Let $n\to\infty$, then $M\to\infty$ and $\eps\downarrow0$. This proves
\begin{equation}\label{eq:transfer-V}
 \max_{j\le k}\left|\frac{L_{I_{n,j}}}{\log n}-\frac A\alpha\right|
 \pto0.
\end{equation}

When $v=\infty$, condition instead on all the depths $H_1,\ldots,H_n$. The selected indices are then fixed, including the tie rule, and independence of the episode pairs allows \eqref{eq:recovery-lower} to be applied to each selected duration. For every $h_0>0$,
\[
 \Pp\bigl(\exists j\le k:
 L_{I_{n,j}}<(\mu^{-1}-\eps)H_{I_{n,j}}\bigr)
 \le\Pp(H_{I_{n,k}}<h_0)+k\sup_{x\ge h_0}r_\eps(x).
\]
Letting $n\to\infty$ and then $h_0\to\infty$ makes the right-hand side vanish. Since $A=1/\mu$, the depth limits and \eqref{eq:transfer-upper} now prove \eqref{eq:transfer-V} in this case too. The marginal maxima and \eqref{eq:transfer-V} give \eqref{eq:light-max} and the duration-ratio limits. Equation~\eqref{eq:transfer-upper} also gives $\limsup_n L_{I_n}/L_n^*\le\rho$ almost surely. An almost-surely convergent subsequence supplied by the probability limit gives equality.

To obtain the duration ranks, choose $\zeta>0$ such that $0<\beta(A/\alpha-\zeta)<\beta(A/\alpha+\zeta)<1$. By \eqref{eq:transfer-V}, with probability tending to one,
\[
 (A/\alpha-\zeta)\log n<L_{I_{n,j}}<(A/\alpha+\zeta)\log n,
 \qquad j\le k.
\]
On this event, the ranks satisfy
\[
 1+\sum_{i\le n}\ind_{\{L_i>(A/\alpha+\zeta)\log n\}}
 \le R_{n,j}
 \le1+\sum_{i\le n}\ind_{\{L_i\ge(A/\alpha-\zeta)\log n\}}.
\]
Equation~\eqref{eq:counts} bounds their logarithmic exponents by $1-\beta(A/\alpha+\zeta)$ and $1-\beta(A/\alpha-\zeta)$. Letting $\zeta\downarrow0$ gives the rank limits in \eqref{eq:rank-limit}. The upper duration bound holds eventually almost surely by \eqref{eq:transfer-upper}, so the lower rank bound also gives
\[
 \liminf_{n\to\infty}\min_{j\le k}\frac{\log R_{n,j}}{\log n}
 \ge1-\rho>0\quad\text{almost surely}.
\]
Consequently the two leading lists are eventually disjoint almost surely.
\end{proof}

\begin{proof}[Proof of Corollary~\ref{cor:overlap} under \eqref{ass:light}]
Choose $\eps>0$ with $A+\eps<\alpha/\beta$, and let $c_\eps>0$ be as in \eqref{eq:joint-long}. Choose $\eta\in(0,1)$ so that $(\alpha+c_\eps)(1-\eta)>\alpha$, and set
\[
 t_n=\log(n/k_n),\qquad
 h_n=\frac{1-\eta}{\alpha}t_n,\qquad
 \ell_n=\frac{1-\eta}{\beta}t_n.
\]
Since $t_n\to\infty$, the marginal tail rates give
\[
 \frac n{k_n}\Pp(H_1\ge h_n)=e^{\eta t_n+o(t_n)},\qquad
 \frac n{k_n}\Pp(L_1>\ell_n)=e^{\eta t_n+o(t_n)}.
\]
The two marginal counts are binomial with means much larger than $k_n$. Chernoff's lower-tail bound shows that both counts are at least $k_n$ except on an event of probability at most $2e^{-ck_n}$, for some $c>0$ and all sufficiently large $n$. On the complementary event, every episode in the overlap satisfies $H_i\ge h_n$ and $L_i>\ell_n$. Since $\ell_n>(A+\eps)h_n$, \eqref{eq:joint-long} gives
\[
 \frac n{k_n}\Pp(H_1\ge h_n,L_1>\ell_n)
 \le C_\eps\exp\!\left[
  \left(1-\frac{(\alpha+c_\eps)(1-\eta)}{\alpha}\right)t_n\right]
 \longrightarrow0.
\]
The joint count is therefore binomial with mean $o(k_n)$. Applying Chernoff's upper-tail bound to this count, we obtain, for each fixed $\delta\in(0,1)$ and all sufficiently large $n$,
\[
 \Pp(O_{n,k_n}>\delta)\le C_\delta e^{-c_\delta k_n},
\]
where $c_\delta,C_\delta>0$. This proves \eqref{eq:overlap}. If $k_n/\log n\to\infty$, the bound is summable in $n$; Borel--Cantelli, applied to $\delta=1/j$, $j\ge2$, gives $O_{n,k_n}\asto0$.
\end{proof}

\section{Proofs for heavy-tailed losses}\label{sec:heavy}

Throughout this section, assume \eqref{eq:common} and \eqref{ass:heavy}. We first obtain a weighted bound on depth--duration deviations and use it to derive a pathwise approximation over the sample. Comparing this error with the gaps between the leading depths identifies the same indices in the same order with probability tending to one. The bound also controls joint exceedance counts, which determine the overlap of growing lists.

\subsection{Depth--duration deviations}

The height and length of the excursion in \eqref{eq:aux} satisfy
\begin{equation}\label{eq:known-heavy}
 \Pp(H^0>x)\sim\Ee\tau\,\overline F(x),\qquad
 \Pp(\tau>t)\sim\Ee\tau\,\overline F(\mu t).
\end{equation}
The stopped-maximum and first-passage results are due to Foss and Zachary~\cite{FossZachary2003} and Denisov and Shneer~\cite{DenisovShneer2013}, respectively.

Conditioning \eqref{eq:known-heavy} on $Y_1>0$ and using \eqref{eq:mean-cycle} gives
\begin{equation}\label{eq:actual-heavy}
 \Pp(H_1>x)\sim\ccyc\overline F(x),\qquad
 \Pp(\mu L_1>x)\sim\ccyc\overline F(x).
\end{equation}
To obtain an almost-sure approximation over the sample, we must control episodes whose depth--duration discrepancy is comparable to the preceding maximum. Summing these probabilities leads to the reciprocal-tail weight in the next lemma.

\begin{lemma}\label{lem:heavy-deviations}
Under \eqref{eq:common} and \eqref{ass:heavy}, put $B=\max(H_1,\mu L_1)$. For each $\eps>0$,
\begin{equation}\label{eq:cone-moment}
 \Ee\!\left[
 \frac{\ind_{\{|H_1-\mu L_1|>\eps B\}}}{\overline F(B)}
 \right]<\infty.
\end{equation}
\end{lemma}
\begin{proof}
Use $(H^0,\tau)$ and put $B=\max(H^0,\mu\tau)$. Fix $\eps\in(0,1)$ and $p\in(1,\min(2,\gamma))$. Regular variation and \eqref{eq:known-heavy} give $\Ee[(Y_1^+)^p]<\infty$ and $\Ee\tau^p<\infty$. Choose $\eta,\delta>0$ with $2\eta+\mu\delta<\eps/2$, and set
\[
 m_h=\lceil2h/\mu\rceil,\qquad
 g(h)=\Pp\!\left(\max_{j\le m_h}|Z_j+\mu j|>\eta h\right).
\]
Truncate the increments at $|Y_i|\le h$ and center them. The omitted mean contributes $o(h)$ over $m_h$ steps, so a union bound and Kolmogorov's maximal inequality give, for large $h$,
\[
 g(h)\le Ch\Pp(|Y_1|>h)
       +\frac C h\Ee[Y_1^2;\ |Y_1|\le h].
\]
By Tonelli's theorem,
\[
 \sum_{m\ge0}\left(2^m\Pp(|Y_1|>2^m)
       +2^{-m}\Ee[Y_1^2;\ |Y_1|\le2^m]\right)
 \le C\Ee|Y_1|<\infty,
\]
and consequently $\sum_m g(2^m)<\infty$.

To control the path before its first large loss, choose $r$ with $r\eta>\gamma+1$ and $r/3>\gamma+1$, and then $\kappa>0$ with $r\kappa<p-1$. Put $\theta=r\log h/h$ and $J=\inf\{j\ge1:Y_j>\kappa h\}$. Using $e^u-1-u\le C_pu^pe^u$ for $u\ge0$ and dominated convergence on $Y_1\le0$ gives
\[
 \begin{aligned}
 \Ee[e^{\theta Y_1};\ Y_1\le\kappa h]
 &\le1-\mu\theta+o(\theta)
     +C\theta^pe^{\theta\kappa h}\Ee[(Y_1^+)^p]\\
 &\le1-\mu\theta/2,
 \end{aligned}
\]
since $\theta^{p-1}e^{\theta\kappa h}\to0$. The process $e^{\theta Z_n}\ind_{\{J>n\}}$ is therefore a nonnegative supermartingale. Its maximal inequality and the bound at time $\lfloor h/\mu\rfloor$ yield
\[
 \begin{aligned}
 \Pp(\max_{n<J}Z_n>\eta h)&\le h^{-r\eta},\\
 \Pp(\tau>h/\mu,\ J>h/\mu)
 &\le(1-\mu\theta/2)^{\lfloor h/\mu\rfloor}
 \le h^{-r/3}.
 \end{aligned}
\]

On $h<B\le2h$ we have $\tau\le m_h$. If $J>\tau$, then either $H^0>h$ or $\tau>h/\mu$; the last two bounds cover both possibilities. The first large loss occurs after $\delta h$ but before recovery with probability at most
\[
 \Pp(\delta h<J<\tau)
 \le\overline F(\kappa h)\sum_{j>\delta h}\Pp(\tau\ge j)
 \le C\overline F(\kappa h)h^{1-p}.
\]
Likewise, independence of the future increments gives
\[
 \begin{aligned}
 &\Pp\!\left(\exists j\le\tau:Y_j>\kappa h,\
       \max_{s\le m_h}|Z_{j+s}-Z_j+\mu s|>\eta h\right)\\
 &\qquad\le\Ee\tau\,\overline F(\kappa h)g(h).
 \end{aligned}
\]

Outside these events, a path with $h<B\le2h$ has its first large loss at a time $j\le\delta h$, and its height before that loss is at most $\eta h$. Put $z=Z_j>0$ and $t=\tau-j$. The future-path bound, together with $Z_\tau\le0<Z_{\tau-1}$, gives
\[
 0\le H^0-z\le\eta h,\qquad
 z-\mu t\le\eta h,\qquad z-\mu(t-1)>-\eta h.
\]
Thus
\[
 |H^0-\mu\tau|
 \le |H^0-z|+|z-\mu t|+\mu j
 \le2\eta h+\mu+\mu\delta h<\eps h
\]
for large $h$. Since $\overline F(\kappa h)\le C\overline F(h)$, the preceding bounds imply
\[
 \Pp(h<B\le2h,\ |H^0-\mu\tau|>\eps h)
 \le C\overline F(h)\bigl(g(h)+h^{1-p}\bigr)
       +Ch^{-r\eta}+Ch^{-r/3}.
\]
Regular variation gives $\overline F(2h)\ge c\overline F(h)$ for large $h$. Hence, for a sufficiently large integer $m_0$,
\[
 \begin{aligned}
 \Ee\!\left[\frac{\ind_{\{|H^0-\mu\tau|>\eps B\}}}{\overline F(B)}\right]
 &\le C+C\sum_{m\ge m_0}
 \frac{\Pp(2^m<B\le2^{m+1},\ |H^0-\mu\tau|>\eps2^m)}
      {\overline F(2^m)}\\
 &<\infty.
 \end{aligned}
\]
Regular variation gives $1/\overline F(h)\le Ch^{\gamma+1}$ for large $h$. Together with $\sum_mg(2^m)<\infty$, $p>1$, and $r\eta,r/3>\gamma+1$, this makes the series converge. Conditioning on $Y_1>0$ proves \eqref{eq:cone-moment}.
\end{proof}

\subsection{Agreement of the leading ranks}

\begin{proof}[Proof of Theorem~\ref{thm:heavy}]
Write $\overline G(x)=\Pp(H_1\ge x)$, put $B_i=\max(H_i,\mu L_i)$, and fix $\eps\in(0,1)$. Call episode $i$ exceptional if $|H_i-\mu L_i|>\eps B_i$. Independence from the preceding episodes and Tonelli's theorem give, with $H_0^*=0$,
\[
 \begin{aligned}
 &\sum_{i\ge1}\Pp\bigl(i\text{ is exceptional},\ B_i>\eps H_{i-1}^*\bigr)\\
 &=\Ee\!\left[\ind_{\{|H_1-\mu L_1|>\eps B_1\}}
       \sum_{i\ge1}(1-\overline G(B_1/\eps))^{i-1}\right]\\
 &=\Ee\!\left[
 \frac{\ind_{\{|H_1-\mu L_1|>\eps B_1\}}}{\overline G(B_1/\eps)}
 \right]<\infty.
 \end{aligned}
\]
The last inequality follows from $\overline G(x)\sim\ccyc\overline F(x)$, regular variation and Lemma~\ref{lem:heavy-deviations}; bounded values of $B_1$ contribute a finite amount. Borel--Cantelli therefore gives, for all sufficiently large $i$, almost surely,
\[
 |H_i-\mu L_i|\le
 \begin{cases}
   \eps H_{i-1}^*,&i\text{ is exceptional},\\
   \dfrac{\eps}{1-\eps}H_i,&i\text{ is not exceptional}.
 \end{cases}
\]
For the second bound, use $|H_i-\mu L_i|\le\eps B_i$ and $H_i\ge(1-\eps)B_i$. Since $H_n^*\to\infty$ almost surely, the finitely many initial terms vanish after normalization, giving
\[
 \limsup_{n\to\infty}
 \frac{\max_{i\le n}|H_i-\mu L_i|}{H_n^*}
 \le\frac{\eps}{1-\eps}\quad\text{almost surely}.
\]
Let $\eps\downarrow0$ through a countable sequence to obtain \eqref{eq:uniform-alignment}. In particular,
\[
 \frac{\mu L_n^*}{H_n^*}\asto1,\qquad
 \frac{\mu L_{I_n}}{H_n^*}\asto1,
\]
which proves the almost-sure duration-ratio limit in \eqref{eq:heavy-conseq}.

To identify the leading indices, we compare the approximation error with the gaps between the largest depths. Equation~\eqref{eq:actual-heavy} gives $n\Pp(H_1>b_nx)\to x^{-\gamma}$ for each $x>0$. For $s$ disjoint intervals in $(0,\infty)$ bounded away from zero, let $p_{n,j}$ be the probability that $H_1/b_n$ lies in interval $j$ and write $np_{n,j}\to\lambda_j$. The multinomial generating function satisfies
\[
 \left(1+\sum_{j=1}^s p_{n,j}(z_j-1)\right)^n
 \longrightarrow\exp\!\left(\sum_{j=1}^s\lambda_j(z_j-1)\right).
\]
Thus the counts converge jointly to independent Poisson variables with mean measure $\gamma x^{-\gamma-1}\,dx$. For positive thresholds $x_1,\ldots,x_{k+1}$,
\[
 \begin{aligned}
 &\Pp(H_{I_{n,j}}\le b_nx_j\text{ for }j\le k+1)\\
 &\qquad=\Pp\!\left(
   \sum_{i\le n}\ind_{\{H_i>b_nx_j\}}\le j-1\text{ for }j\le k+1\right).
 \end{aligned}
\]
The count limits therefore determine the joint order-statistic law. Under $x\mapsto x^{-\gamma}$, the limiting mean measure becomes Lebesgue measure on $(0,\infty)$, giving
\begin{equation}\label{eq:height-orders}
 (H_{I_{n,j}}/b_n)_{j=1}^{k+1}
 \dto(\Gamma_j^{-1/\gamma})_{j=1}^{k+1}.
\end{equation}

Equation~\eqref{eq:height-orders} makes $H_n^*/b_n$ tight. Together with \eqref{eq:uniform-alignment}, this gives
\begin{equation}\label{eq:sample-error}
 \frac1{b_n}\max_{i\le n}|H_i-\mu L_i|\pto0.
\end{equation}
The limiting upper order statistics in \eqref{eq:height-orders} are strictly positive and strictly decreasing. Define the smallest leading gap
\[
 \Delta_{n,k}=\min_{1\le j\le k}
      (H_{I_{n,j}}-H_{I_{n,j+1}}).
\]
Then
\begin{equation}\label{eq:gap-small}
 \lim_{\delta\downarrow0}\limsup_{n\to\infty}
 \Pp(\Delta_{n,k}\le2\delta b_n)=0.
\end{equation}
On the event
\[
 \max_{i\le n}|H_i-\mu L_i|<\delta b_n,
 \qquad \Delta_{n,k}>2\delta b_n,
\]
for every $j\le k$ and every index $i$ outside the first $j$ depth positions,
\[
 \mu L_{I_{n,j}}-\mu L_i
 \ge H_{I_{n,j}}-H_i-2\delta b_n
 \ge\Delta_{n,k}-2\delta b_n>0.
\]
Hence the first $k$ ordered index lists coincide. Equations \eqref{eq:sample-error}--\eqref{eq:gap-small}, followed by $\delta\downarrow0$, prove \eqref{eq:topk}.

The uniform error in \eqref{eq:sample-error} and \eqref{eq:height-orders} give \eqref{eq:order-limit} by Slutsky's theorem. Taking $k=1$ in \eqref{eq:topk} gives $\Pp(R_n=1)\to1$.
\end{proof}

\begin{proof}[Proof of Corollary~\ref{cor:overlap} under \eqref{ass:heavy}]
For growing lists, we use the deviation estimate to compare joint and marginal exceedance probabilities. Put $B=\max(H_1,\mu L_1)$. For each $\eps>0$, Lemma~\ref{lem:heavy-deviations} gives
\[
 \begin{aligned}
 &\frac{\Pp(B>x,\ |H_1-\mu L_1|>\eps B)}{\overline F(x)}\\
 &\quad\le\Ee\!\left[
   \frac{\ind_{\{|H_1-\mu L_1|>\eps B\}}}{\overline F(B)};\ B>x\right]
 \longrightarrow0.
 \end{aligned}
\]
For $\eta>0$, the event $H_1>(1+\eta)x$, $\mu L_1\le x$ implies $B=H_1$ and $|H_1-\mu L_1|>\eta B/(1+\eta)$. Hence
\[
 \begin{aligned}
 \Pp(H_1>x,\mu L_1>x)
 &\ge\Pp(H_1>(1+\eta)x)\\
 &\quad-\Pp(B>x,\ |H_1-\mu L_1|>\eta B/(1+\eta)).
 \end{aligned}
\]
Together with \eqref{eq:actual-heavy} and regular variation, this gives
\[
 (1+\eta)^{-\gamma}
 \le\liminf_{x\to\infty}
 \frac{\Pp(H_1>x,\mu L_1>x)}{\ccyc\overline F(x)}
 \le\limsup_{x\to\infty}
 \frac{\Pp(H_1>x,\mu L_1>x)}{\ccyc\overline F(x)}\le1.
\]
Letting $\eta\downarrow0$ proves $\Pp(H_1>x,\mu L_1>x)\sim\ccyc\overline F(x)$.

Choose $u_n\to\infty$ with $n\ccyc\overline F(u_n)\sim k_n$, using a generalized inverse of the regularly varying tail. For $\eps\in(0,1)$, choose $\delta>0$ so that $1-\eps<(1+\delta)^{-\gamma}<1$, and put $x_n=(1+\delta)u_n$. The two marginal counts above $x_n$ and the joint count are each binomial, with mean asymptotic to $k_n(1+\delta)^{-\gamma}$. Chernoff's bounds show that both marginal counts are at most $k_n$ and the joint count is at least $(1-\eps)k_n$, except on an event of probability at most $C_\eps e^{-c_\eps k_n}$ for some $c_\eps,C_\eps>0$ and all sufficiently large $n$. Whenever both marginal counts are at most $k_n$, every joint exceedance belongs to both leading lists, so
\[
 O_{n,k_n}\ge\frac1{k_n}\sum_{i\le n}
       \ind_{\{H_i>x_n,\ \mu L_i>x_n\}}.
\]
Consequently,
\[
 \Pp(O_{n,k_n}<1-\eps)\le C_\eps e^{-c_\eps k_n}.
\]
This proves \eqref{eq:overlap}. If $k_n/\log n\to\infty$, summability and Borel--Cantelli, applied to $\eps=1/j$, $j\ge2$, give $O_{n,k_n}\asto1$.
\end{proof}

\section{Deterministic observation horizons}\label{sec:calendar}

When investment performance is assessed over a prescribed horizon of $N$ periods, the number of completed drawdowns is random. Both ranking laws continue to hold as $N\to\infty$. Define
\begin{equation*}
 K_N=\max\{i\ge0:r_i\le N\},\qquad N\in\Nn.
\end{equation*}
Since $r_i=\sum_{j=1}^i(W_j+L_j)$, the renewal strong law gives
\begin{equation}\label{eq:renewal}
 \frac{K_N}{N}\asto\frac1{\ccyc}.
\end{equation}
Since $K_N$ depends on the durations being ranked, we use deterministic bounds on $K_N$ to transfer the probability limits.\footnote{When $K_N<k$, missing indices and associated statistics are assigned fixed values, with missing ranks set to one and ratios with no completed episode set to zero. Set $O_{K_N,k}=0$ when $K_N<k$.}

\begin{corollary}[Observation-time rank laws]\label{cor:calendar}
Under \eqref{eq:common} and \eqref{ass:light},
\begin{equation*}
 \begin{gathered}
 \frac{L_{I_{K_N}}}{L^*_{K_N}}\pto\rho,
 \qquad \frac{\log R_{K_N}}{\log N}\pto1-\rho,\\[5pt]
 \frac{H^*_{K_N}}{\log N}\asto\frac1\alpha,
 \qquad \frac{L^*_{K_N}}{\log N}\asto\frac1\beta,
 \qquad \frac{L_{I_{K_N}}}{\log N}\pto\frac A\alpha.
 \end{gathered}
\end{equation*}
The limits for the deepest episode's duration and rank hold jointly for every fixed number of deepest episodes. For every fixed $k$, the index sets of the $k$ deepest and $k$ longest episodes are eventually disjoint almost surely. Under \eqref{eq:common} and \eqref{ass:heavy}, for every fixed integer $k\ge1$,
\begin{equation}\label{eq:calendar-alignment}
 \Pp(I_{K_N,j}=J_{K_N,j}\ \text{for all }1\le j\le k)
 \longrightarrow1.
\end{equation}
The duration ratio satisfies
\[
 \frac{L_{I_{K_N}}}{L^*_{K_N}}\asto1.
\]
Choose $a_N\to\infty$ so that $N\overline F(a_N)\to1$. Then
\begin{equation}\label{eq:calendar-heavy}
 \left(\frac{H_{I_{K_N,j}}}{a_N},
       \frac{\mu L_{I_{K_N,j}}}{a_N}\right)_{j=1}^{k}
 \dto\left(\Gamma_j^{-1/\gamma},\Gamma_j^{-1/\gamma}\right)_{j=1}^{k}.
\end{equation}
For every deterministic $k_N\to\infty$ with $k_N=o(N)$, the overlap $O_{K_N,k_N}$ converges in probability to zero under \eqref{ass:light} and to one under \eqref{ass:heavy}.
\end{corollary}
The factor $\ccyc$ in the depth tail of a completed nonzero episode cancels with the completion rate $1/\ccyc$. The observation-time scale is therefore the negative-increment quantile determined by $N\Pp(X_1<-a_N)\to1$. For the shifted Pareto model in Section~\ref{sec:examples}, one may take $a_N=N^{1/\gamma}-c$ for all sufficiently large $N$, so $a_N\sim N^{1/\gamma}$.

\begin{proof}
Fix $\delta\in(0,1)$, and put
\[
 n_- =\left\lfloor\frac{(1-\delta)N}{\ccyc}\right\rfloor,
 \qquad
 n_+ =\left\lceil\frac{(1+\delta)N}{\ccyc}\right\rceil.
\]
By \eqref{eq:renewal}, $\Pp(n_-\le K_N\le n_+)\to1$. The first $k$ depth indices remain unchanged throughout $n_-\le m\le n_+$ unless a new observation enters the leading list. For $i>n_-$, exchangeability with independent random tie breaking gives probability $k/i$ of entering the first $k$ positions. Breaking ties by the earlier index can only decrease this probability for the new observation. Hence
\[
 \Pp\bigl(\exists m\in[n_-,n_+],\ j\le k:
       I_{m,j}\ne I_{n_-,j}\bigr)
 \le k\sum_{i=n_-+1}^{n_+}\frac1i
 \longrightarrow k\log\frac{1+\delta}{1-\delta}.
\]
The same bound holds for the duration indices $J_{m,j}$.

In the light-tailed case, \eqref{eq:transfer-V} at the deterministic sample size $n_-$ and $\log n_-/\log N\to1$ give, for every $\eps>0$,
\[
 \limsup_{N\to\infty}\Pp\!\left(
   \max_{j\le k}\left|\frac{L_{I_{K_N,j}}}{\log N}-\frac A\alpha\right|
   >\eps\right)
 \le k\log\frac{1+\delta}{1-\delta}.
\]
Letting $\delta\downarrow0$ proves the selected duration limits. The marginal maximum limits and eventual separation transfer pathwise from Theorem~\ref{thm:light}, since $K_N\to\infty$ and $\log K_N/\log N\to1$ almost surely. The duration-ratio limits follow. The count limits in \eqref{eq:counts} also hold pathwise along $K_N$; bracketing the selected durations as in Section~\ref{sec:light-ranks} therefore gives the rank limits.

In the heavy-tailed case, the duration-ratio limit follows pathwise from Theorem~\ref{thm:heavy} and $K_N\to\infty$. The two leading lists agree at $n_-$ with probability tending to one. If neither list changes between $n_-$ and $n_+$, they also agree at $K_N$ on the renewal bracketing event. The entry bounds therefore give
\[
 \limsup_{N\to\infty}
 \Pp\bigl(\exists j\le k:I_{K_N,j}\ne J_{K_N,j}\bigr)
 \le2k\log\frac{1+\delta}{1-\delta}.
\]
Letting $\delta\downarrow0$ proves \eqref{eq:calendar-alignment}.

For \eqref{eq:calendar-heavy}, the tail normalization is
\[
 n_\pm\ccyc\overline F(a_Nx)
 \longrightarrow(1\pm\delta)x^{-\gamma},\qquad x>0.
\]
The same Poisson-count argument as in the proof of Theorem~\ref{thm:heavy} shows that the first $k$ depth values at $n_\pm$, divided by $a_N$, converge to
\[
 \bigl((1\pm\delta)^{1/\gamma}\Gamma_j^{-1/\gamma}\bigr)_{j=1}^k.
\]
Each ordered depth is nondecreasing in sample size. On the bracketing event, the ordered vector at $K_N$ is therefore bounded coordinatewise by those at $n_-$ and $n_+$. Taking their joint distribution functions and then letting $\delta\downarrow0$ proves
\[
 (H_{I_{K_N,j}}/a_N)_{j=1}^k
 \dto(\Gamma_j^{-1/\gamma})_{j=1}^k.
\]
Finally, \eqref{eq:uniform-alignment} and tightness of $H_{n_+}^*/a_N$ give
\[
 \frac1{a_N}\max_{i\le n_+}|H_i-\mu L_i|
 =\frac{\max_{i\le n_+}|H_i-\mu L_i|}{H_{n_+}^*}
   \frac{H_{n_+}^*}{a_N}\pto0.
\]
This comparison applies to all selected indices on the bracketing event, giving the full joint vector limit \eqref{eq:calendar-heavy}.

For the overlap, let $k_N\to\infty$ with $k_N=o(N)$. Eventually $k_N\le n_-$. A replacement in either leading list changes the overlap proportion by at most $1/k_N$. Combining this with the entry bound $k_N/i$ for each list gives
\[
 \Ee\max_{n_-\le m\le n_+}
     |O_{m,k_N}-O_{n_-,k_N}|
 \le2\sum_{i=n_-+1}^{n_+}\frac1i
 \longrightarrow2\log\frac{1+\delta}{1-\delta}.
\]
Corollary~\ref{cor:overlap} applies at the deterministic size $n_-$ because $k_N=o(n_-)$. Combining its limit with the renewal bracket and Markov's inequality, then letting $\delta\downarrow0$, proves the overlap limits.
\end{proof}

\section{Conclusion}

For light-tailed losses with a finite exponential moment, the deepest drawdown is eventually different from the longest almost surely. Reaching a large depth and sustaining a long episode have different exponential costs, and the deepest episode's duration rank tends to infinity. Under regular variation, the deepest drawdown is also the longest with probability tending to one. A common large loss determines depth and first-order recovery time: the duration ratio tends to one almost surely, and the leading rankings agree in order with probability tending to one for every fixed list size.

This contrast persists when the list size grows but its fraction of the sample vanishes: the overlap proportion tends in probability to zero in the light-tailed case and to one in the heavy-tailed case. Both ranking laws apply to completed episodes sampled by count or within a deterministic observation horizon.

\bibliographystyle{elsarticle-num-noids}
\bibliography{references}
\end{document}